\documentclass[10pt]{article}

\usepackage{amsthm,amsmath,amssymb}
\usepackage{enumerate}

\def\lra{\longrightarrow}

\def\ol{\overline}

\def\rb\}{\right\}}

\def\wt{\widetilde}
\def\x{\times}
\def\hook{\hfill\hbox{\rlap{$-$}\;\raise .6ex\hbox{$\lrcorner$}}}

\def\del{\partial}

\renewcommand{\phi}{\varphi}

\renewcommand{\epsilon}{\varepsilon}
\def\eps{\varepsilon}

\newtheorem{thm}{Theorem}[section]
\newtheorem{prop}[thm]{Proposition}
\newtheorem{lem}[thm]{Lemma}
\newtheorem{cor}[thm]{Corollary}

\newtheorem{rem}[thm]{Remark}
\newtheorem{dEf}[thm]{Definition}

\newtheorem{exs}[thm]{Examples}
\newtheorem{question}{Question}

\newcommand{\comment}[1]{}

\def\supp{\text{{\rm supp}}\,}

\def\iNt{\text{{\rm int}}}

\def\bbn{{\mathbb N}}

\def\bbq{{\mathbb Q}}
\def\bbr{{\mathbb R}}

\def\bbz{{\mathbb Z}}

\def\ca{{\mathcal A}}

\def\cc{{\mathcal C}}
\def\cd{{\mathcal D}}

\def\cf{{\mathcal F}}

\def\co{{\mathcal O}}

\def\cs{{\mathcal S}}
\def\ct{{\mathcal T}}

\title{Relative structure cycles and the
existence of smooth Lyapunov 1-forms for flows}
\author{Janko Latschev\footnote{Postal address:
    Institut f\"ur Mathematik, Humboldt Universit\"at zu Berlin, Unter
    den Linden 6, D-10099 Berlin, Germany, \tt
    latschev@mathematik.hu-berlin.de}}  
\begin{document}
\maketitle

%
%

\section{Introduction}

In this work we address the following 

\begin{question}
\label{quest:dual}
Given a smooth vector field $V$ on a closed manifold $M$, what are
necessary and sufficient conditions for the existence of a smooth {\em
closed} 1-form $\lambda$ on $M$ with the same zero set as $V$, such
that outside this common zero set $\lambda(V)< 0$? 
\end{question}

Obviously, the cohomology class of such a form $\lambda$ evaluates
negatively on all homology classes of non-constant periodic orbits of
$V$. In particular, if $\lambda$ is exact, then there can be no
non-constant periodic orbits. The main result below asserts that a
suitable generalization of this necessary condition on the cohomology
class of $\lambda$ is also sufficient.  

The history of this question goes back to Schwartzman's beautiful work
\cite{Sw}. One of his results answers the above question in the case
that $V$ has no zeros. Schwartzman associates 
to every finite positive invariant Borel measure on $M$ a
1-dimensional homology class which he calls {\em asymptotic cycle} of
the measure. Then he proves that there exists a closed form $\lambda$
with $\lambda(V)<0$ on $M$ in a given cohomology class $\xi\in
H^1(M;\bbz)$ if and only if $\xi$ evaluates negatively on all these
asymptotic cycles. This was one of the inputs for Sullivan's
influential paper \cite{Su}, where he treats existence of closed forms
transverse to quite general cone structures, including for example
foliations of all dimensions.
In \cite{Fa:1,Fa:2}, Farber studied the dynamics of gradient-like
vector fields for closed 1-forms with isolated zeros and introduced
the term Lyapunov 1-form. This concept was developed further in 
joint works with Farber, Kappeler and Zehnder \cite{FKLZ:1,FKLZ:2},
where we arrived at the following definition. It generalizes the
notion of a smooth Lyapunov function in the sense of Conley \cite{Co:1}. 
\begin{dEf}
\label{def:lyap}
Given a vector field $V$ on a manifold $M$ generating a global
flow  and a flow-invariant closed subset $Z\subset M$, we say that a
smooth, closed 1-form $\lambda$ is a {\em Lyapunov 1-form for $(V,Z)$}
if
\begin{enumerate}[(L1)]
\item $\lambda(V)< 0$ on $M \setminus Z$, and
\item $\lambda$ vanishes pointwise on $Z$ and is exact in a
neighborhood of $Z$.
\end{enumerate}
\end{dEf}
It was shown in \cite{FKLZ:2} that the exactness assumption in (L2) is
redundant if either $Z$ is a Euclidean neighborhood retract or the
cohomology class of $\lambda$ is rational. In terms of this
definition, our question can be rephrased as
\begin{question}
\label{quest:lyap}
Given a smooth vector field $V$ on the closed manifold $M$, a closed,
flow-invariant subset $Z\subset M$ and a cohomology class $\xi\in
H^1(M;\bbr)$, what are necessary and sufficient conditions for the
existence of a Lyapunov 1-form $\lambda$  for $(V,Z)$ representing
$\xi$?
\end{question}
While this may appear to be more general than Question~\ref{quest:dual}, 
it is clear that the existence of a Lyapunov 1-form depends only on
the oriented foliation of $M \setminus Z$ by the flow lines of $V$ and
the way that $Z$ is embedded in $M$. Thus multiplying the vector field
of Question~\ref{quest:lyap} by a nonnegative function whose zero set
coincides with $Z$, we return to the previous situation.

Various consequences of the existence of Lyapunov 1-forms for the flow
generated by $V$ are known. Generalizing earlier work of
Farber~\cite{Fa:1}, in \cite{La:1} the author proved that if the set
$Z$ is small (meaning that its Lusternik-Schnirelman category is
smaller than Farber's category associated to the cohomology class
$\xi$ of a Lyapunov 1-form), there must exist a cyclically ordered
chain of orbits outside $Z$ such that the forward limit set of one is
contained in the same connected component of $Z$ as the backward limit
set of the next. Also, under suitable conditions on $Z$, Fan and
Jost~\cite{FJ:1} established generalized Morse-Novikov inequalities, 
relating the topology of $Z$ to the Novikov numbers associated to the
cohomology class of any Lyapunov 1-form for $(V,Z)$. 

The main goal of this work is to give answers to
Question~\ref{quest:lyap}. For that purpose, we introduce the notion
of {\em measures coherent relative to $Z$}. These are flow-invariant,
positive, locally finite measures $\mu$ on $M \setminus Z$, satisfying
the additional condition that 
\begin{equation*}
\int_{M \setminus Z} df(V)\,d\mu =0
\end{equation*}
whenever $f$ is a smooth function with $df\equiv 0$ on some
neighborhood of $Z$.
Now suppose that $Z$ is an isolated invariant set in the sense of
Conley \cite{Co:1}, i.e. a closed invariant subset which is the
maximal set with this property in some closed neighborhood of itself.   
In this case, 
coherent measures have an alternative description, which
roughly asserts that for each connected component of $Z$ the amount of
mass ``flowing in'' equals the amount of mass ``flowing out'' (see
\S~\ref{flowcycles} for the precise statement and some
examples). Denote by $H_Z\subset H^1(M;\bbr)$ the subspace 
of cohomology classes vanishing in some neighborhood of $Z$. In other
words, $H_Z$ is the kernel of the restriction map from $H^1(M;\bbr)$
to the \v Cech cohomology $\check H^1(Z;\bbr)$ of $Z$. Then it is easy
to see that every coherent measure $\mu$ gives rise to a homomorphism
$A_\mu:H_Z \lra \bbr$, defined as
\begin{equation}
A_{\mu} (\xi) := \int_{M \setminus Z} \alpha(V) d\mu,
\end{equation}
where $\alpha$ is any representative of $\xi$ vanishing near $Z$. As
this homomorphism equals Schwartzman's asymptotic cycle of the
measure when $Z$ is empty, we call it {\em relative
asymptotic cycle of $\mu$}. With these concepts, our first main
theorem reads as follows. 
\begin{thm}
\label{thm:lyap_exist_general}
Let $V$ be a smooth vector field on the closed manifold $M$ and let $Z$
be an isolated invariant set for the flow of $V$. Then there exists a
Lyaponov 1-form for $(V,Z)$ representing the cohomology class $\xi\in
H_Z\subset H^1(M;\bbr)$ if and only if the relative asymptotic cycle
of every coherent measure $\mu$ satisfies  
\begin{equation}
\label{eq:1}
A_\mu(\xi) < 0. 
\end{equation}
\end{thm}
In \S~\ref{flowforms}, we derive the following consequence of this
theorem, which completely answers Question~\ref{quest:dual} in case
the zero set $Z$ of $V$ is finite.
\begin{cor}
\label{cor:finite_z}
Let $V$ be a smooth vector field on the closed manifold $M$ with
finite zero set $Z=\{z_1,\dots,z_l\}$. Then there exists a Lyapunov
1-form for $(V,Z)$ if and only if each $z_i$ is an isolated invariant
set and the relative asymptotic cycle $A_\mu:H^1(M;\bbr) \lra \bbr$ of
every coherent measure $\mu$ satisfies (\ref{eq:1}).
\end{cor}
The proof of Theorem~\ref{thm:lyap_exist_general} follows the same
basic scheme as pursued by Sullivan \cite{Su}, which consists of two
steps. In a first step, he uses the Hahn-Banach theorem to relate the
existence question for transversal closed forms to the homological
properties of structure cycles of the considered cone structure. These
are closed currents, i.e. elements of the dual space of the space of
smooth differential forms vanishing on all exact forms. In the second
step, Sullivan identifies these structure cycles in specific
situations. For example, in the case of foliations, structure cycles
are in bijective correspondence to transverse invariant measures (see
also \cite{RS:1}). 

In \S~\ref{cone-struct} we extend the basic definitions and
the first technical theorem of \cite{Su} to our relative situation. 
Here the main difficulty is to deal with the non-compactness of the
set $M \setminus Z$ where the inequality $\lambda(V)<0$ is supposed to
hold. However, we show that with the appropriate setup a suitable
Hahn-Banach argument can still be applied. Theorem~\ref{thm:main1} of
this section is stated for cone structures in the generality considered
by Sullivan, which are however defined on the difference of two compact
sets. We expect this result to have applications beyond the case of
flows considered here.

After recalling some structural results from \cite{Su} in
\S~\ref{struct_fol}, in \S~\ref{flowcycles} we consider flows and
prove (see Proposition~\ref{prop:flow_cycles}) that if $Z$ is an 
isolated invariant set, then our relative structure cycles correspond
to (suitable restrictions of) coherent measures. Together with
Theorem~\ref{thm:main1} this result immediately yields
Theorem~\ref{thm:lyap_exist_general}. 

Note that given Lyapunov 1-forms $\lambda_1$ and $\lambda_2$ for
$(V,Z_1)$ and $(V,Z_2)$ respectively, the form $\lambda_1+\lambda_2$
clearly is a Lyapunov 1-form for $(V,Z_1 \cap Z_2)$. 
For integral cohomology classes, any representative gives rise to a
map $M \lra S^1$. Using these observations, we prove in
\S~\ref{flowforms} the following necessary and sufficient condition
for the existence of a Lyapunov 1-form in a rational cohomology class
for an arbitrary closed invariant subset $Z\subset M$.  
\begin{thm}
\label{thm:lyap_exist_integral}
Let $V$ be a smooth vector field on a closed manifold $M$ and let
$Z\subset M$ be a closed subset invariant under its flow. Then there
exists a Lyaponov 1-form for $(V,Z)$ representing 
$\xi\in H^1(M;\bbq)$ if and only if $Z$ can be written as a
countable intersection of isolated invariant sets $Z = \cap_{i\in I}
Z_i$ such that for each $i\in I$ there exists a Lyapunov 1-form for
$(V,Z_i)$ representing $\xi$. 
\end{thm}
Now we consider the following 
\begin{question}
\label{quest:smallest}
If Lyapunov one-forms representing a given cohomology class $\xi\in
H^1(M;\bbr)$ exist for a given vector field $V$ with respect to some
$Z$, then what is the smallest subset $Z_\xi(V)\subset M$ relative to 
which such a form can be found?
\end{question}
Before stating our answer in the case of rational cohomology classes,
we introduce a few definitions.
Recall that, given a flow $\phi$ on $M$, an $(\eps,T)$-chain for
$\phi$ from $x\in M$ to $y\in M$ is a finite sequence
$x=x_1,\dots,x_{n+1}=y$ of points in $M$ and real numbers
$t_1,\dots,t_n\geq T$ such that the distance $d(x_i\cdot
t_i,x_{i+1})<\eps$ for all $i=1,\dots,n$. This distance is 
measured with respect to some auxilliary Riemannian metric. 
The {\em chain recurrent set $R=R(\phi)$} is the set of all points
$x\in M$ such that for every $\eps>0$ and $T>1$ there is an
$(\eps,T)$-chain from $x$ to itself. It is a closed, flow invariant
set that contains all nonwandering points. Note that if $\eps$ is
sufficiently small (e.g. smaller than the injectivity radius), then
one can associate to every $(\eps,T)$-chain from $x$ to itself a
singular $1$-cycle, which is given as the sum for $i=1,\dots,n$ of the
oriented flow lines from $x_i$ to $x_i\cdot t_i$ followed by the
unique shortest curves joining $x_i\cdot t_i$ to $x_{i+1}$. 
In \cite{FKLZ:1}, we introduced the following notion.
\begin{dEf}
\label{def:r_xi}
Let a flow $\phi$ on $M$ and $\xi\in H^1(M;\bbr)$ be given. The subset
$R_\xi(\phi)\subset R(\phi)$ of the chain recurrent set consists of
those points $x\in M$ such that for all sufficiently small $\eps>0$
and all $T>1$ there exists an $(\eps,T)$-chain from $x$ to itself such
that the homology class $z\in H_1(M;\bbz)$ of the associated cycle satisfies 
$\langle \xi,z \rangle =0$.
\end{dEf}
Alternatively, $R_\xi(\phi)$ can be characterized as the set of points
$x\in R(\phi)$ whose lifts to the abelian cover $M_\xi$ of $M$
associated to $\xi\in H^1(M;\bbr)$ are chain recurrent for the lifted
flow. This fact, along with some examples, is discussed in
\cite{FKLZ:1}. In \S \ref{flowforms}, we prove the following result.
\begin{thm}
\label{thm:characterize_R_xi}
Let $V$ be a smooth vector field with flow $\phi$ on the closed
manifold $M$ and let the class $\xi\in H^1(M; \bbq)$ be given.
If the collection $\{Z_i\}$ of all closed $\phi$-invariant sets such
that there exists a smooth Lyapunov 1-form for $(V,Z_i)$ representing the
class $\xi$ is nonempty, then 
\begin{equation}
\bigcap_{i\in I} Z_i = R_\xi(\phi).
\end{equation}
\end{thm}
We note in passing that a fundamental theorem of Conley guarantees
existence of an exact Lyapunov 1-form for $(V,R(\phi))$. The standard
proof of this fact uses a representation of $R(\phi)$ as above, where the
$Z_i$ run through all attractor-repeller pairs for the flow in $M$
(see \cite{Co:1} and also \cite{FKLZ:2}). 

In \cite{FKLZ:2}, we established necessary and sufficient conditions
for the existence of a smooth Lyapunov 1-form for $(V,R_\xi(\phi))$ in
a given cohomology class $\xi\in H^1(M;\bbr)$, assuming that
$R(\phi)\setminus R_\xi(\phi)$ is closed. One of several equivalent
conditions there was given in terms of Schwartzman's asymptotic
cycles. In \S~\ref{flowforms}, we show how this result can be
recovered in the present framework. More precisely, we reprove
\begin{prop}\cite{FKLZ:2}
\label{prop:lyap_exist_r_xi}
Let $V$ be a smooth vector field on the smooth closed manifold $M$ and
let $\xi\in H^1(M;\bbr)$ be given. Denote by $R_\xi$ the part of the
chain recurrent set $R$ introduced in Definition~\ref{def:r_xi} and
assume that $R\setminus R_\xi$ is closed and that $\xi$ vanishes on
some neighborhood of $R_\xi$.  

Then there exists a Lyapunov 1-form for $(V,R_\xi)$ representing $\xi$
if and only if the asymptotic cycle $A_\mu$ of every finite positive
invariant measure $\mu$ whose support intersects $M \setminus R_\xi$
satisfies 
\begin{equation*}
A_\mu(\xi) < 0.
\end{equation*}
\end{prop}

{\bf Acknowledgements: } It is a pleasure to thank B. Lawson for
pointing me towards \cite{Su}, as well as K.~Cieliebak,
M.~Farber, T.~Kappeler and E.~Zehnder for stimulating
discussions. During the finishing stage of
this work, the author held a position financed by the
DFG-Schwerpunktprogramm ``Differentialgeometrie''.

%
%


%
%

\section{Structure currents for relative cone structures}
\label{cone-struct}

Let $M$ be a smooth manifold. In this section, we modify the basic
definitions of \cite{Su} for cone structures defined on differences
$X\setminus Z$ of compact subsets of $M$ and extend Sullivan's basic
existence theorem for transversal closed forms to this more general context.  
Throughout, we denote by $\cd^p$ the space of compactly supported
smooth $p$-forms with the $C^\infty$-topology, and by $\cd_p$ the
space of $p$-dimensional currents. With their natural topologies,
these are dual topological vector spaces \cite{dR}. 

\begin{dEf}
A $p$-dimensional cone structure $\cs_Z$ on the compact subset $X\subset M$
relative to the closed subset $Z\subset X$ consists of a compact convex cone
$C_p(x) \subset  \Lambda_pT_xM$ (based at $0$) for each $x\in X\setminus Z$
which depends continuously on $x$.  
\end{dEf}
Recall that a convex cone $C \subset H$ based at $0$ in a topological vector
space is called compact if there exists a continuous linear functional $L:H
\lra \bbr$ such that the set $L^{-1}(1) \cap C$ is compact in $H$. 
Continuity is measured with respect to Hausdorff convergence of the set of
unit $p$-vectors in $C_p$ in some auxilliary chart near $x\in M\setminus Z$.  

\begin{dEf}
Given a $p$-dimensional relative cone structure $\cs_Z$ on $X \subset
M$, a smooth $p$-form $\lambda\in \cd^p$ satisfying 
\begin{enumerate}[(i)]
\item $\lambda(v)>0$ whenever $x\in X\setminus Z$ and $v\in C_p(x)$, and
\item $\lambda = 0$ pointwise on $Z$ 
\end{enumerate}
is called {\em transversal relative to $Z$} or simply {\em relatively
tranversal}.  
\end{dEf}
We use this sign convention here to be consistent with \cite{Su}. Note
the following existence statement: 
\begin{prop}(cf. \cite[I.4]{Su})
\label{prop:exist1}
Any relative cone structure $\cs_Z$ on a compact subset $X\subset M$
admits a relatively transversal form. 
\end{prop}

\begin{proof} 
Fix a sequence of smooth functions $f_n:M \lra [0,1]$ such that $Z_{1/2n}
\subset f_n^{-1}(0) \subset Z_{1/n}$, where $Z_\eps$ denotes the set of points
whose distance to $Z$ is smaller than $\eps$. Arguing as in the proof of
\cite[I.4]{Su}, for each $n\in \bbn$ we construct a smooth form
$\omega_n$ which is transversal for the restriction of $\cs_Z$ to
the compact complement of $Z_{1/2n}$, which is a cone structure in
Sullivan's sense. Then, for a sufficiently rapidly decaying sequence
$c_n>0$ of real numbers, the form
$$
\omega = \sum_{n= 1}^\infty c_n \cdot f_n \omega_n
$$
is smooth and transversal relative to $Z$. 
\end{proof}
In our definition of relative structure currents and relative structure
cycles we will use an approach essentially dual to that of
Sullivan. Given a relative cone structure $\cs_Z$, for each open
neighborhood $U$ of $Z$ we introduce the open cone 
$$
\co_U := \{ \alpha \in \cd^p(M) \; | \; \alpha(v) > 0 \text{ \rm for
all } v \in \cs_{Z|M \setminus U} \}. 
$$
Moreover, for each cohomology class $\xi\in H^1(M;\bbr)$ we introduce
the closed cone 
\begin{eqnarray*}
\ca_{U,\xi} := &\{ \alpha \in \cd^p(M) \; | \; \alpha \equiv 0 \text{
\rm on } Z, \alpha(v) \geq 0 \text{ \rm for all } v \in \cs_{Z|\ol U},\\
 & d\alpha=0, [\alpha]=\xi\,\}. 
\end{eqnarray*}
We have seen in Proposition \ref{prop:exist1} above that $\co_U$ is
non-empty for every $U$. $\ca_{U,0}$ is also non-empty, as it contains
the linear subspace 
$\cf_U := \{ d\beta \,|\, d\beta \equiv 0 \text{ \rm on } \ol U \}$.
Similarly, if $\xi$ vanishes in a neighborhood of $\ol U$, the cone
$\ca_{U,\xi}$ is non-empty because there exists a representative
$\alpha_0$ of $\xi$ vanishing on $\ol U$. In fact in this case we have
$\ca_{U,\xi}= \alpha_0 + \ca_{U,0}$. 
The following easy observation will turn out to be quite useful.

\begin{lem}
\label{lem:lyap_exist}
There exists a relatively transversal closed form  representing the
class $\xi\in H^p(M; \bbr)$ if and only if  
$$
\co_U \cap \ca_{U,\xi}\neq \varnothing
$$ 
for all neighborhoods $U$ of $Z$.
\end{lem}

\begin{proof}
The necessity of the condition is immediate from the definitions,
because any relatively transversal closed form representing $\xi$ is
in the intersection of both cones.

Conversely, fix a countable system of neighborhoods $U_i$ of $Z$ with
$\cap U_i = Z$ and choose $\alpha_i \in \co_{U_i} \cap
\ca_{U_i,\xi}$. Then for a sufficiently rapidly decaying sequence of
constants $c_i$, the sum 
$$
\alpha := \sum c_i \alpha_i
$$
is a well-defined smooth closed form. From the definitions we see that it
satisfies both properties required of relatively transversal forms. Moreover,
the cohomology class of $\alpha$ equals $c\cdot \xi$, where $c = \sum
c_i$. So setting $\lambda := \frac 1 c \alpha$, we obtain the required
transversal form representing $\xi$. 
\end{proof}
\begin{rem}
\label{rem:small_nbhd}
Note that if $W\subset U$ are two neighborhoods of $Z$, then
\begin{equation*}
\co_W \cap \ca_{W,\xi} \subset \co_U \cap \ca_{U,\xi}
\end{equation*}
for all cohomology classes $\xi\in H^p(M;\bbr)$. In particular, if
$\co_U \cap \ca_{U,\xi} = \varnothing$, then the same is true for all
smaller neighborhoods $W$.
\end{rem}
\begin{dEf}
Given an open neighborhood $U$ of $Z$, a current $T$ is called a
{\em structure current relative to $U$} if $T(\alpha)>0$ for all
$\alpha\in \co_U$. 
If in addition $T(\alpha) \leq 0$ for all $\alpha\in \ca_{U,0}$, then
$T$ is called a {\em structure cycle relative to $U$}.
\end{dEf}
By a relative stucture current (resp. cycle) we simply mean a
structure current (resp. cycle) relative to some neighborhood $U$ of $Z$.

Denote by $H_U\subset H^p(M;\bbr)$ the linear subspace of cohomology
classes vanishing on some neighborhood of $\ol U$. Then every
structure cycle $T$ relative to $U$ induces a homomorphism $\ct: H_U
\lra \bbr$ given by
$$
\ct(\xi) = T(\alpha),
$$
where $\alpha$ is any representative of $\xi$ vanishing on some
neighborhood of $\ol U$. As the difference of two such representatives
is an element of the linear subspace $\cf_U \subset \ca_{U,0}$ on
which $T$ vanishes, this is indeed well-defined.
Note that if $V \subset U$ is another neighborhood of $Z$, then
clearly $H_U \subset H_V$. Thus, if $H^p(M;\bbr)$ is finite-dimensional
(e.g. if $M$ is compact), then for all sufficiently small
neighborhoods $U$ of $Z$ the subspaces $H_U$ coincide with the
subspace $H_Z\subset H^p(M;\bbr)$ of classes vanishing in some
neighborhood of $Z$.
The following statement is our generalization of Theorem I.7 in
\cite{Su} to the relative case. 

\begin{thm}
\label{thm:main1}
Let $M$ be a smooth manifold and let $\cs_Z$ be a $p$-dimensional cone
structure on the compact subset $X\subset \iNt(M)$ relative to the
compact subset $Z\subset X$. 

Then a cohomology class $\xi\in H_Z \subset H^p(M;\bbr)$
contains a relatively transversal form if and only if for every open
neighborhood $U$ of $Z$ with $\xi\in H_U$ and every structure cycle
$T$ relative to $U$ we have 
\begin{equation}
\ct(\xi) >0.
\end{equation}
\end{thm}

\begin{proof}
If there exists no relatively transversal closed form representing the
class $\xi$, then according to Remark \ref{rem:small_nbhd} for all
sufficiently small neighborhoods $U$ the cones $\co_U$ and
$\ca_{U,\xi}$ are disjoint. If $\xi\in H_Z$, we can assume $U$ to be so
small that $\xi$ vanishes on some neighborhood of $\ol U$, so that
$\ca_{U,\xi}\neq \varnothing$. By the
Hahn-Banach theorem, there exists a continuous linear functional
$T:\cd^p \lra \bbr$ with $T_{\co_U}>0$ and $T_{|\ca_{U,\xi}}\leq 0$. Any
such $T\in \cd_p$ is clearly a structure current. Moreover, given
$\alpha\in \ca_{U,0}$, we have $\beta +c \cdot \alpha\in \ca_{U,\xi}$
for every $\beta \in \ca_{U,\xi}$ and $c>0$. As $T(\beta)$ is finite
and $T(\beta +c \cdot \alpha)\leq 0$, this forces $T(\alpha)\leq 0$,
so that $T$  is in fact a structure cycle relative to $U$ with
$\ct(\xi) \leq 0$.  

Conversely, assume there exists a closed relatively
transversal form $\lambda$ representing a class $\xi\in H_Z$.
Let $U$ be an open neighborhood of $Z$ with  $\xi\in H_U$ and let $T$
be a structure cycle relative to $U$. Given any representative
$\alpha'$ of $\xi$ vanishing on some neighborhoohd of $\ol U$, we
compute  
$$
\ct(\xi) = T(\alpha') = T(\lambda) - T(\lambda-\alpha') >0
$$
because $\lambda\in \co_U$ and $\lambda-\alpha'\in \ca_{U,0}$. This
completes the proof.
\end{proof}
We explicitly remark a special case.
\begin{cor}
\label{cor:exact}
Under the assumptions of Theorem~\ref{thm:main1}, an exact relatively
transversal form exists if and only if there are no relative structure
cycles. $\Box$
\end{cor}

%
%

\section{Local study of relative structure cycles}
\label{struct_fol}

We begin this section by showing that our definition of
relative structure currents reduces to the one of Sullivan 
\cite{Su} in the absolute case $Z=\varnothing$.
\begin{lem} 
\label{lem:struct}
Let $U$ be some open neighborhood of $Z$.
\begin{enumerate}[(i)]
\item The cone of structure currents relative to $U$ coincides with
  the complement of $0$ in the closed convex cone $\cc_U\subset \cd_p$
  generated by evaluations on elements $v\in \cs_{Z|M\setminus U}$. 
\item For structure cycles relative to $U$ we have 
$$
\supp(\del T) \subset \del U = \ol U \cap (M \setminus U). 
$$
\end{enumerate}
\end{lem}

\begin{proof}
We first prove part (i). Indeed, if there existed
some structure current $T \notin \cc_U$, then by the Hahn-Banach
theorem there would exist some smooth one-form $\alpha\in \cd^p(M)$
such that $\alpha_{|\cc_U} \geq 0$ and $a=\alpha(T) <0$. Let $\beta$
be some element of $\co_U$, and let $b=T(\beta)>0$. Then the 
form $\gamma= \alpha - \frac a {2b} \beta$ satisfies $\gamma_{|\cc_U}
> 0$ and $\gamma(T)= a - \frac a 2<0$. As $\cc_U$ contains evaluation
on every $p$-vector $v\in \cs_{Z|M\setminus U}$, we see that $\gamma \in 
\co_U$. But this contradicts the fact that $T_{|\co_U}> 0$, so that
the assumption that a structure current $T\notin \cc_U$ existed is
disproven. The converse assertion is immediate.

Part (ii) follows directly from the fact that $\ca_{U,0}$ contains the
linear subspace $\cf_U$ introduced above, on which all structure
cycles vanish. Indeed, any $(p-1)$-form $\beta$ vanishing on some
neighborhood of $\del U$ can be decomposed as a sum of two such forms
$\beta_1$ and $\beta_2$ with $\supp \beta_1 \subset U$ and $\supp
\beta_2 \subset M \setminus \ol U$. As $\supp T \subset M \setminus
U$, we have $T(d\beta_1)=0$. Since $d\beta_2\in \cf_U$, we also have
$T(d\beta_2)=0$. This proves that $\supp \del T$ is contained in every
neighborhood of $\del U$, proving (ii). 
\end{proof}

From Lemma~\ref{lem:struct}~(i) and the Riesz representation theorem, we
immediately deduce (cf. Proposition~I.8 in \cite{Su})

\begin{prop}
\label{prop:rep_measure}
Any structure current $T$ relative to $U$ can be represented as
\begin{equation}
T(\alpha) = \int_{X \setminus U} \alpha(v)\, d\mu_T,
\end{equation}
where $\mu_T$ is a non-negative measure on $X\setminus U$ and $v$ is a
$\mu$-integrable map into $p$-vectors satisfying $v(x)\in C_p(x)$
$\mu$-a.e. $\square$
\end{prop}

Next we record, for later use, Sullivan's local description 
for structure cycles in case that the cone structure is determined
by an oriented foliation of a compact manifold. In fact, the proof of
Theorem~I.12 of \cite{Su} applies directly to establish the
following result.

\begin{prop}
\label{prop:cycle_local}
Suppose the relative cone structure $\cs_Z$ on the compact manifold
$M$ is determined by an oriented foliation of $M\setminus Z$ of class $C^1$.
Fix a neighborhood $U$ of $Z$ and a closed neighborhood $B\subset M
\setminus \ol U$ foliated by closed disks $L_y$ on leaves of the 
foliation, and denote by $D$ the space of these disks.

Then for every foliation current $T$ relative to $U$ with $\supp T
\cap B \neq \varnothing$ and $\supp \del T \cap B = \varnothing$
there is a non-negative measure $\nu$ on $D$ (which is
unique in the interior of $D$), so that on the interior of $B$ the
current $T$ can be represented as 
\begin{equation}
T(\alpha) = \int_D [L_d](\alpha)\,d\nu(d),
\end{equation}
where $[L_d]$ denotes the current of integration over the disk $L_d$
oriented by the foliation. $\square$
\end{prop}
We conclude this section with some examples of relative structure
cycles.
\begin{exs}
\begin{enumerate}[(i)]
\item \cite{Su} For a relative cone structure $\cs_Z$ induced form a
foliation of $M \setminus Z$, any transversal invariant measure with
compact support in $M \setminus Z$ determines a relative structure cycle.
\item For flows, the simplest example of a relative structure cycle
that is not the restriction of an absolute structure cycle is the
current of integration over $\gamma \cap (M \setminus U)$, where
$\gamma$ is a homoclinic orbit, i.e. $\omega(\gamma) = \alpha(\gamma)
= p \in Z$ and $U$ is some neighborhood of $Z$.
\item Consider the foliation of $\bbr^n$ by $k$-dimensional subspaces
parallel to $\bbr^k\x \{0\}$. This foliation can be pulled back to
$S^n$ via stereographic projection map $p:S^n \setminus\{p\} \lra
\bbr^n$ from some point $p\in S^n$ to give a foliation of $S^n
\setminus \{p\}$. Denote by $\cs$ the associated cone structure
relative to $Z=\{p\}$. Then given any neighborhood $U \supset Z$, the
intersection of any leaf with $S^n \setminus U$ is a structure cycle
relative to $U$. In this example Corollary~\ref{cor:exact} asserts the
elementary fact that if $1\leq k < n$, then no closed form transversal
relative to $Z=\{p\}$ exists.
\end{enumerate}
\end{exs}

%
%

\section{Structure cycles for 1-dimensional singular foliations}
\label{flowcycles}

In this section, we discuss structure cycles for oriented
1-dimensional singular foliations. Let $V$ be a smooth vector field on
a compact manifold $M$ and suppose $Z\subset M$ is a closed subset invariant
under the flow of $V$ which contains all zeros of $V$. We consider the
relative cone structure $\cs_Z$ determined by the restriction of $V$ to
$M \setminus Z$ and begin with an elementary result on supports of
structure cycles. 
\begin{lem}
\label{lem:support}
Let $V$ be a  smooth vector field on the compact manifold $M$ with flow
$\phi$ and suppose $Z$ is contained in the chain recurrent set $R(\phi)$. 
Then the support of every relative structure cycle is also contained
in $R(\phi)$. 
\end{lem}
\begin{proof}
Denote by $L:M \lra \bbr$ a smooth Lyapunov function in the sense of
Conley for the flow $\phi$. This means that $dL$ vanishes pointwise on
the chain recurrent set $R(\phi)$ and $dL(V) < 0$ on $M\setminus
R(\phi)$. In particular, if $Z \subset R(\phi)$, then $\alpha:=-dL\in
\ca_{U,0}$ for every neighborhood $U$ of $Z$. Now if $T$ is any
structure current relative to $U$ whose support intersects $M
\setminus R(\phi)$, then $T(\alpha)>0$ and so $T$ cannot be a relative
structure cycle. This proves the lemma.
\end{proof}
We next introduce a class of measures which will turn out to be
closely related to relative structure cycles.
\begin{dEf}
We say that a positive, locally finite Borel measure $\mu$ on $M
\setminus Z$ is {\em coherent relative to $Z$} if 
\begin{enumerate}
\item[(C)]\label{eq:coherent}
for all smooth functions $f:M \lra \bbr$ with $df\equiv 0$ in some
neighborhood of $Z$ we have
$$
\int_{M \setminus Z} df(V)\,d\mu = 0.
$$
\end{enumerate}
\end{dEf}
In the discussion so far, the open sets $U$ were essentially 
arbitrary neighborhoods of $Z$. We will now see
that more can be said about the structure of coherent measures (and
structure cycles relative to $U$) in the special case that $\ol U$ is
an isolating block. 

Recall that a smooth hypersurface $S\subset M$ is called a {\em local cross
section} for the flow if for some $\delta>0$ the flow induced map $S \x
(-\delta,\delta) \lra M$ is a diffeomorphism onto the open set $S\cdot
(-\delta,\delta)$.  

\begin{dEf} (cf. \cite{Ch})
The closure $B$ of an open set in $M$ is called an {\em isolating block for
the set $Z$} with respect to the flow $\phi$ if there exist local cross
sections $S_\pm \subset M$ and $\delta>0$ as above such that
\begin{enumerate}[(B1)]
\item the sets $\del_\pm B :=S_\pm \cap B$ are closed in $M$.
\item $S_+ \cdot (-\delta,\delta) \cap B = \del_+B \cdot [0,\delta)$ and
$S_- \cdot (-\delta,\delta) \cap B = \del_-B \cdot (-\delta,0]$, 
\item for every  $x\in \del B \setminus (\del_+B \cup
\del_-B)$ there are $t_1< 0$ and $t_2> 0$ such that $x\cdot [t_1,t_2]
\subset \del B$ and $x\cdot t_1\in \del_+B$, $x\cdot t_2 \in \del_-B$, and
\item $Z$ is the maximal invariant subset of $B$.
\end{enumerate}
\end{dEf}

Note that the maximal invariant subset $Z$ of an isolating block $B$ is
always contained in the interior of $B$, as is immediate from (B2) and
(B3). Thus it is an isolated invariant set in the sense of Conley
\cite{Co:1}. Conversely, it is known that every neighborhood of an isolated
invariant set contains an isolating block \cite{Ch}. 

Let $B$ be an isolating block for $Z$ and let $\mu$ be a positive Borel
measure on $M \setminus Z$ which is locally finite and invariant under the
flow. By the standard argument, $\mu$ induces slice measures $\nu_\pm$
on $S_\pm$ and hence, by restriction, on $\del_\pm B$
(cf.~Proposition~\ref{prop:cycle_local}). Moreover, we observe the
following fact.
\begin{lem}
\label{lem:boundary}
Given an invariant, positive, locally finite Borel measure $\mu$ on $M
\setminus Z$ and an isolating block $B$ for $Z$, we have
\begin{equation}
\int_{M \setminus B} df(V)\,d\mu 
= \int_{\del_+ B} f d\nu_+ - \int_{\del_- B} f d\nu_-,
\end{equation}
where $\nu_\pm$ are the slice measures on $\del_\pm B$ induced by the measure
$\mu$.
\end{lem}
\begin{proof}
We compute
\begin{eqnarray*}
\int_{M \setminus B} df (V) \, d\mu 
&=& \int_{M \setminus B} \lim_{t\to 0} \frac 1 t \big(
f(x\cdot t) -f(x)\big)\, d\mu(x) \\
&=& \lim_{t\to 0} \Bigg[ \frac 1 t \int_{ (M \setminus B) \cap (M \setminus
B)\cdot t} f(x\cdot t) -f(x) \, d\mu(x) \\
& &  + \frac 1 t \int_{\del_+B \cdot [-t,0]} f(x\cdot t) \,d\mu(x) 
- \frac 1 t \int_{\del_-B \cdot [0,t]} f(x \cdot t) \,d\mu(x) \Bigg].\\
\end{eqnarray*}
The first summand in this expression vanishes because the measure
$\mu$ is invariant under the flow. As $f$ is smooth, the sum of the
last two terms converges to   
$$
\int_{\del_+B} f(b) \,d\nu_+(b) - \int_{\del_-B} f(b)\,d\nu_-(b).
$$
This completes the proof of Lemma \ref{lem:boundary}.
\end{proof}
Inside $\del_\pm B$, we have the subset $a_\pm \subset \del_\pm B$
consisting of those points whose forward (resp.~backward) limit set is
contained in $B$ and thus in $Z$. We denote the complements $\del_\pm B
\setminus a_\pm$ by $\del^*_\pm B$, respectively. It is well known that the
flow-induced map $h_B: \del^*_+ B \lra \del^*_-B$ is a
homeomorphism. Moreover, from the invariance of the measure $\mu$ we
deduce that the map $h_B: (\del^*_+B,\nu_+) \lra (\del^*_-B,\nu_-)$ is
measure-preserving. 

Given a connected component $Z_0\subset Z$, we denote by $a_\pm(Z_0) \subset
\del_\pm B$ the set of points whose forward (resp. backward) limit set lies
in $Z_0$. With this notation, we obtain the following equivalent
characterization of coherent measures.
\begin{lem}
\label{lem:coherent_equiv}
A locally finite, positive, invariant Borel measure $\mu$ is coherent
relative to an isolated invariant set $Z$ if and only if 
\begin{enumerate}
\item[(C')]\label{eq:coherent2}
for every isolating block $B$ and every connected component $Z_0$ of
$Z$ we have 
\begin{equation*}
\nu_+(a_+(Z_0)) = \nu_-(a_-(Z_0)).
\end{equation*}
\end{enumerate}
\end{lem}
\begin{rem}
Note that if $B$ and $B_1$ are two isolating blocks for $Z$ with
$B \subset \iNt(B_1)$, by the invariance of the measure $\mu$ condition
$(C')$ is true for $B$ if and only if it is true for $B_1$. 
Thus to check coherence it suffices to check condition $(C')$ in some
isolating block for $Z$.
\end{rem}
\begin{proof}
First suppose $\mu$ satisfies condition $(C')$ 
for every isolating block $B$ and every component $Z_0\subset Z$. Given a
smooth function $f:M \lra \bbr$ with $df\equiv 0$ on the neighborhood $W$
of $Z$, choose an isolating block $B\subset W$. Then 
\begin{eqnarray*}
\int_{M\setminus Z} df(V)\,d\mu &=& \int_{M \setminus B} df(V)\,d\mu\\
&=& \int_{\del_+B} f\, d\nu_+ - \int_{\del_-B} f \,d\nu_-\\
&=& 0,
\end{eqnarray*}
because $\mu$ is invariant and the function $f$ takes the same value
on corresponding points of $\del_\pm B$. 

Conversely, assume that $\mu$ satisfies condition $(C)$. 
Given some isolating block $B$ for $Z$, for every component $Z'\subset
Z$ we introduce the sets
\begin{eqnarray*}
A_+(Z') &:=& \{ x\in B \,|\, x\cdot [0,\infty) \subset B \text{ \rm and
  } \ol{x \cdot [0,\infty)} \setminus x \cdot[0,\infty) \subset Z' \}, \\
A_-(Z') &:=& \{ x\in B \,|\, x\cdot (-\infty,0] \subset B \text{ \rm and
  } \ol{x \cdot(-\infty,0]} \setminus x \cdot (-\infty,0] \subset Z'
  \},  \\ 
A(Z') &:=& A_+(Z') \cup A_-(Z').
\end{eqnarray*}
The following assertion is proven in the appendix. 
\begin{lem}
\label{lem:approx}
Given any connected component $Z_0$ of $Z$, the characteristic
function of $A(Z_0)$ is the pointwise limit of smooth functions $f_n$
with $df_n(V)\equiv 0$ on $B$ and $df_n\equiv 0$ on some neighborhood
of $Z$ (which depends on $n$).  
\end{lem}
Using this sequence of functions we compute
\begin{eqnarray*}
\nu_+(a_+(Z_0))-\nu_-(a_-(Z_0)) 
&=& \int_{\del_+B} \chi_{A(Z_0)} \,d\nu_+ 
- \int_{\del_-B} \chi_{A(Z_0)} \,d\nu_- \\
&=& \lim_{n\to\infty} \int_{\del_+B} f_n \,d\nu_+ 
- \int_{\del_-B} f_n \,d\nu_-\\
&=& \lim_{n\to\infty} \int_{M} df_n(V)\,d\mu = 0.
\end{eqnarray*}
This shows (C') and thus completes the proof of Lemma~\ref{lem:coherent_equiv}.
\end{proof}
\begin{exs}
\begin{enumerate}[(i)]
\item A finite invariant measure with support disjoint from $Z$
is coherent relative to $Z$.
\item Let $\gamma_1,\dots,\gamma_k,\gamma_{k+1}=\gamma_1$ be a chain
of orbits outside $Z$ such that the forward limit set of $\gamma_i$ is
contained in the same connected component of $Z$ as the backward limit
set of $\gamma_{i+1}$ for all $1\leq i\leq k$. Each $\gamma_i$ carries a
natural invariant measure $\mu_i$, namely the pushforward of Lebesgue measure
on $\bbr$ by the map $t \mapsto x_i\cdot t$ for any point $x_i\in
\gamma_i$. The sum $\mu=\sum_{i=1}^k \mu_i$ is a coherent measure
relative to $Z$.
\item Suppose the flow of the nowhere vanishing vector field $V_1$ has
$M_1\subset M$ as a minimal set. Let $\mu_1$ be an invariant measure with
$\supp \mu =M_1$. Define a new vector field $V= f \cdot V_1$, where
$f:M \lra [0,1]$ is smooth with zero set $Z\not\supset M_1$. Then $\mu
:= \frac 1 f \cdot \mu_1$ is a coherent measure relative to $Z$. 
\end{enumerate}
\end{exs}
As mentioned in the introduction, every coherent measure $\mu$ gives
rise to a homomorphism $A_\mu:H_Z \lra \bbr$ via 
\begin{equation}
\label{eq:def_a_mu}
A_\mu(\xi) := \int_{M \setminus Z} \alpha(V)\; d\mu,
\end{equation}
where $\alpha$ denotes any representative of the class $\xi$
vanishing near $Z$. As the integral on the right vanishes on exact
forms with support in $M \setminus Z$ because the measure is
invariant, $A_\mu$ is indeed well-defined. 
In the case $Z=\varnothing$, every finite positive invariant measure
$\mu$ is coherent, and the corresponding homomorphism
$A_\mu:H^1(M;\bbr)\lra \bbr$ is nothing but Schwartzman's asymptotic
cycle \cite{Sw} associated to the measure $\mu$. For that reason, we
call the homomorphism $A_\mu:H_Z \lra \bbr$ the {\em relative asymptotic
cycle} of the coherent measure $\mu$. 

We now describe the relation between coherent measures relative to
$Z$ and relative structure cycles.
\begin{prop}
\label{prop:flow_cycles}
Let $Z$ be an isolated invariant set for the flow generated by
the vector field $V$ and let $U$ be the interior of an isolating block
for $Z$. 

Then there is a bijective correspondence between coherent measures
$\mu$ relative to $Z$ and structure cycles relative to $U$. 
\end{prop}
\begin{rem}
It will be obvious from the proof that under this correspondence, the
relative asympotic cycle $A_\mu:H_Z \lra \bbr$ of the coherent
measure $\mu$ and the homomorphism $\ct_{U,\mu}$ induced by the
associated structure cycle $T_{U,\mu}$ are related via
$A_{\mu|H_U}=\ct_{U,\mu}$. 
\end{rem}
%
%
%
\begin{proof}
Fix an isolating block $B$ for $Z$ and a coherent measure $\mu$. As
the support of $\mu$ is invariant under the flow, it cannot be
completely contained in $U=\iNt(B)$, because otherwise $Z$ would not be 
the maximal invariant subset in $B$.

We claim that the current $T_{U,\mu}$ defined by
$$
T_{U,\mu} (\alpha) := \int_{M \setminus U} \alpha(V)\, d\mu
$$
is a structure cycle relative to $U$. As the measure $\mu$ is
positive, we clearly have $T_{U,\mu}(\beta)>0$ for all $\beta\in
\co_U$. Moreover, for $df\in \ca_{U,0}$ we compute using
Lemma~\ref{lem:boundary} that
\begin{eqnarray*}
T_{U,\mu}(df) = \int_{M\setminus U} df(V)\,d\mu = \int_{\del_+B} f \,d\nu_+ -
\int_{\del_-B} f\,d\nu_- \leq 0,
\end{eqnarray*}
where we used property (C') of the coherent measure $\mu$, the fact
that $h_B$ is measure preserving and the fact that $f$ is
nondecreasing along flow lines in $B$. 
This completes the proof that $T_{U,\mu}$ is a relative structure cycle.

\bigskip

Conversely, let $T$ be a structure cycle relative to $U$. Recall from
Proposition~\ref{prop:rep_measure} that it is of the form
\begin{equation*}
T(\alpha) = \int_{M \setminus U} \alpha(V)\, d\mu_T,
\end{equation*}
where $\mu_T$ is a non-negative measure on $M\setminus U$. Using that
$T$ vanishes on the subspace $\cf_U\subset \ca_{U,0}$ of differentials
of functions which are locally constant in a neighborhood of $\ol U$,
we find that $\mu_T$ is invariant in the sense that if $W\cdot[0,t]
\subset M \setminus U$ and $W$ is Borel measurable, then
$\mu_T(W)=\mu_T(W\cdot t)$. 
Moreover, using the fact that $\del_\pm B$ are closed subsets of the local
sections $S_\pm$, near which the flow induces a product structure,
Proposition~\ref{prop:cycle_local} gives the existence of measures
$\nu_\pm$ on $\del_\pm B$ such that 
\begin{equation}
T_{|\del_+B\cdot [-\delta,0]}(\alpha) = \int_{\del_+B} \left( 
\int_{b\cdot   [-\delta,0]} \alpha\right) \, d\nu_+(b)
\end{equation}
and a similar expression exists for $T_{|\del_-B \cdot
[0,\delta]}$. Now a computation as in the proof of Lemma 
\ref{lem:boundary} shows that  
\begin{equation}
\del T(f) = \int_{\del_+ B} f d\nu_+ - \int_{\del_- B} f d\nu_-.
\end{equation}
In addition, the vanishing of $T$ on the differentials of functions
which are constant along flow lines inside $U$ and constant on
components of $Z$, combined with Lemma~\ref{lem:approx}, implies
directly that the measure $\mu_T$ has property (C') 
with respect to $B$ and that $h_B:(\del^*_+ B,\nu_+) \lra
(\del^*_-B,\nu_-)$ is measure preserving. But this is all we need to
extend $\mu_T$ invariantly to a measure $\mu$ on all of $M\setminus
Z$. The extension is unique because the difference of two such
extensions would be a (possibly signed) invariant measure with support
in $U \setminus Z$, contradicting the fact that $Z$ is the maximal
invariant subset of the isolating block $B=\ol U$. The measure $\mu$
is locally finite on $M\setminus Z$, because any compact subset of $M
\setminus Z$ can be covered by finitely many sets of the form $(M
\setminus U) \cdot t$ for some $t\in \bbr$ and $\mu_T(M \setminus U)<
\infty$.  

Thus we have associated to every coherent measure a structure cycle
relative to $U$, and to every structure cycle relative to $U$ an
coherent measure. As the two constructions are inverses of
one another, the proof of Proposition \ref{prop:flow_cycles} is complete.
\end{proof}

%
%

\section{Lyapunov 1-forms for flows}
\label{flowforms}

In this section we use the results established so far to prove the
assertions of the introduction. So let $V$ be a $C^1$-vector field,
denote by $\phi$ its flow, let $Z$ be a closed, flow-invariant
subset containing the zeros of $V$. Denote by $\cs_Z$ the relative
cone structure induced by the restriction of $V$ to $M\setminus
Z$.

\subsubsection*{Proof of Theorem~\ref{thm:lyap_exist_general}}
Observe first that $\alpha$ is a closed form transversal relative
to $Z$ in the sense of \S~\ref{cone-struct} whose cohomology class is in
$H_Z$ if and only if $-\alpha$ is a Lyapunov 1-form in the sense of
Definition~\ref{def:lyap}. 

Now suppose that there exists a Lyapunov 1-form $\lambda$ for $(V,Z)$
representing the class $\xi\in H_Z$. Choose an isolating block for $Z$
and denote its interior by $U$. Then in view of
Proposition~\ref{prop:flow_cycles} every coherent measure $\mu$ gives
rise to a structure cycle relative to $U$, and by
Theorem~\ref{thm:main1} the associated relative asymptotic cycle
satisfies $A_\mu(\xi)<0$. 

Conversely, suppose that there is no Lyapunov 1-form for $(V,Z)$
representing $\xi\in H_Z$. Pick an isolating block $B$ for $Z$ such
that the restriction of $\xi$ to $B$ vanishes, and set $U=\iNt(B)$. By
our choice of $B$, both cones $\co_U$ and $\ca_{U,-\xi}$ are
nonempty. In view of Lemma~\ref{lem:lyap_exist} we can assume $U$ to
be so small that these cones are disjoint. Let $T$ be a structure
cycle relative to $U$ separating these cones with $\ct(\xi)\geq 0$,
which exists by
Theorem~\ref{thm:main1}. Proposition~\ref{prop:flow_cycles} and the
remark following now show that this structure cycle corresponds to a
coherent measure $\mu$ with  $A_\mu(\xi)\geq 0$. \hfill $\Box$

\subsubsection*{Proof of Corollary~\ref{cor:finite_z}}

It suffices to prove that if a Lyapunov 1-form $\lambda$ for $(V,Z)$
representing some cohomology class $\xi\in H^1(M;\bbr)$ exists,
then each $z_i\in Z$ is an isolated invariant set. 

So suppose some $z_i\in Z$ is not an isolated invariant set. Pick
a closed ball $B$ centered at $z_i$ and disjoint from $Z \setminus
\{z_i\}$ such that the restriction of the given Lyapunov
1-form $\lambda$ is exact on $B$, i.e. $\lambda_{|B}=dL$. By our
assumption there exists a point $x\in B \setminus \{z_i\}$ whose
orbit is completely contained in $B$. Denote by $X\subset B$ the
closure of the orbit of $x$. Now note that no point $y\in X \setminus
\{z_i\}$ can be the minimum or maximum of the function
$L_{|X}$, because $L(y \cdot (-1))> L(y) > L(y \cdot 1)$. As $X$ is
compact, this gives the desired contradiction. \hfill $\Box$

\comment{
Note that its forward and backward limit sets are closed,
invariant subsets of $B_i$. If one of them does not contain $z_i$,
then it supports a finite, positive, invariant measure $\mu$ for which
$A_\mu(\xi)=0$ because $\lambda$ is exact on $B_i$. In view of
Theorem~\ref{thm:lyap_exist_general}, this is impossible. 
On the other hand, if both forward and backward limit sets of $x$
contain $z_i$, then for each $\eps>0$ we find $t_-(\eps)<0$ and
$t_+(\eps)>1$ such that $d(x \cdot t_\pm,z_i)<\eps$. Now choose $\eps>0$
so small that the integral of $\lambda$ over any path in $B_i$ of
length at most $2\eps$ is smaller than $|\int_{x \cdot [0,1]}
\lambda|$. Choose a path $\gamma'$ of length at most $2\eps$
connecting $x \cdot t_+$ to $x\cdot t_-$. Then
$$
\int_{x\cdot[t_-,t_+]} \lambda + \int_{\gamma'} \lambda <
\int_{x\cdot[0,1]} \lambda + \left| \int_{x\cdot[0,1]} \lambda \right|
=0,
$$
again contradicting the fact that $\lambda$ is exact on $B_i$. This
finishes the proof of Corollary~\ref{cor:finite_z}. $\Box$
}

\subsubsection*{Proofs of Theorems~\ref{thm:lyap_exist_integral} and
\ref{thm:characterize_R_xi}}

We now analyse the special situation when the cohomology class
$\xi$ of a Lyapunov 1-form is integral. We assume without loss of
generality that $M$ is connected, and fix some point $x_0\in M$. Given
a closed, flow-invariant subset $Z\subset M$ and a Lyapunov 1-form
$\lambda$ for $(V,Z)$ representing $\xi\in H^1(M,\bbz)$, the assignment
\begin{eqnarray}
L:M &\lra&  S^1 = \bbr/\bbz\label{eq:map_L}\\
x &\mapsto& \int_{x_0}^x \lambda  \mod 1,\notag
\end{eqnarray}
where the integral is along any smooth path connecting $x_0$ to $x$,
defines a smooth map whose differential is $\lambda$. Using this
simple observation, we will prove Theorem~\ref{thm:lyap_exist_integral} 
and Theorem~\ref{thm:characterize_R_xi}.
\begin{proof} (of Theorem~\ref{thm:lyap_exist_integral})
As the existence of a Lyapunov 1-form depends only on the ray
$\bbq_+\cdot \xi \in H^1(M;\bbq)$, we will assume that the class $\xi$
is integral. 
If $Z$ is the countable intersection of isolated invariant sets for
which Lyapunov 1-forms $\lambda_i$ representing $\xi$ exist, then for a
suitable choice of rapidly decaying, positive coefficients $c_i$
the form $\lambda= \sum c_i\lambda_i$ will be the desired Lyapunov
1-form for $(\phi,Z)$. 

Conversely, given a Lyapunov 1-form $\lambda$ for $(V,Z)$ representing
an integral class $\xi\in H^1(M;\bbz)$, we consider the map $L:M \lra
S^1$ as above. Pick one point $s_i\in S^1$ in each of the at most
countably many connected components of the set of regular values of
$L$. If $(s_i-2\eps_i,s_i+2\eps_i)$ denotes an interval consisting of
regular values, one straightforwardly varifies that the set
$L^{-1}(S^1 \setminus (s_i-\eps,s_i+\eps))$ is an isolating 
block for the flow. Denote by $Z_i$ the corresponding isolated
invariant set. Clearly we have $Z\subset Z_i$ for each
$i$. Conversely, if $x\notin Z$, then the flow line of $x$ crosses at
least one of the hypersurfaces $L^{-1}(s_i)$, and so $x$ cannot be in
the corresponding $Z_i$. This proves $Z= \cap Z_i$. 

It remains to show that for each $i$ there exists a Lyapunov 1-form for
$(V,Z_i)$ representing $\xi$. So let $\mu$ be some measure
which is coherent relative to $Z_i$. As the support of $\mu$ cannot be
contained in the isolating block, it must intersect
$L^{-1}(s_i)$. Let $\wt \alpha=a(t)dt$ be a form
on $S^1$ supported in $(s_i-\eps,s_i+\eps)$, such that $a(t) \geq 0$,
$a(s_i)>0$ and $\int_{S^1}\wt \alpha = 1$.  Then the form
$\alpha=L^*(\wt\alpha)=a(L(x))\cdot \lambda$
represents $\xi$, has support outside the isolating block for $Z_i$
and satisfies $\alpha(V)\leq 0$ there, with $\alpha(V)<0$ on $L^{-1}(s_i)$.
Now it follows directly from the definition (\ref{eq:def_a_mu}) that
$A_\mu(\xi)<0$. Hence the condition of
Theorem~\ref{thm:lyap_exist_general} holds for $Z_i$, and the proof of  
Theorem~\ref{thm:lyap_exist_integral} is complete. 
\end{proof}
Next we prove Theorem~\ref{thm:characterize_R_xi}.
\begin{proof}
(of Theorem~\ref{thm:characterize_R_xi}) We will show that
\begin{enumerate}[(a)]
\item if there exists a Lyapunov 1-form for $(V,Z)$ representing
  $\xi\in H^1(M;\bbz)$, then $R_\xi(\phi) \subset Z$, and
\item if there exists a Lyapunov 1-form for $(V,Z)$ representing $\xi$
  for some closed, $\phi$-invariant subset, then there exists one for
  $(V,R_\xi(\phi))$.
\end{enumerate}
This clearly suffices to establish the Theorem.

Again, given our Lyapunov 1-form for $(V,Z)$ representing $\xi$, we
consider the map $L:M \lra S^1$ as in (\ref{eq:map_L}) above. To prove
(a), we fix a point $x\in M \setminus Z$, and our goal is to show that
$x\notin R_\xi$. Replacing it by some other point on its orbit if
necessary, we may assume that $l=L(x)$ is a regular value of $L$. As
the set of regular values of $L$ is open ($M$ is compact), we can fix
$0<\eps<1/4$ such that $[l-2\eps,l]\subset S^1$ consists only of
regular values. This means that the function $dL(V)=\lambda(V)$ is
strictly negative on $L^{-1}([l-2\eps,l])$, and so it admits a
negative upper bound there. 

In particular, there exists some $T>1$ such that $\int_y^{y\cdot
T}\lambda < - \eps$ for all $y\in L^{-1}([l-\eps,l])$. Furthermore, as
$L$ is uniformly continuous, we can fix $\delta>0$ such that if
$d(y_1,y_2)<\delta$, then the two points can be joined by a unique
geodesic and $|L(y_1)-L(y_2)|< \eps/2$.

Now let $(z_0=x,z_1,\dots,z_n=x,t_1,\dots,t_n)$ be any finite
$(\delta,T)$-chain from $x$ to itself. As $z_0=x$, we see from our
choice of $T$ that $\int_{z_0}^{z_0\cdot t_1} \lambda<-\eps$, so that
$\int_{z_0}^{z_1} \lambda <-\eps/2$. In particular, either
$d(z_1,z_0)>\delta$ or the integral along the path obtained from
following the orbit of $z_0$ until time $t_1$ and then connecting to
$z_1$ by the unique geodesic from $z_0\cdot t_1$ is $\leq
-1+\eps/2$. Continuing this argument inductively, we see that the 
integral of $\lambda$ along the closed curve associated to the
$(\delta,T)$-chain is $\leq -1+\eps/2$ (in fact $-1$, because $\xi$ is an
integral class). As the $(\delta,T)$-chain from $x$ to itself was
arbitrary, we conclude that $x$ cannot be in $R_\xi$. This proves (a).

To prove (b), we split $M$ open along the inverse image $L^{-1}(l)$ of
some regular value of $L$ to obtain a manifold $N$ with two boundary
components and an obvious smooth projection $p:N \lra M$. As $V$ is
transverse to $L^{-1}(l)$, it points inward at one boundary component
of $N$, which we denote by $N_1$, and outward at the other component,
which we call $N_0$. Pick a smooth function $\rho:N \lra [0,1]$ with
$\rho^{-1}(0) =N_0\cup N_1$ and consider the smooth vector field $V_N
:= \rho\cdot p^{-1}_*(V)$. It induces a complete flow $\Phi$ on $N$
for which $N_1$ is a repeller and $N_0$ is an attractor in the sense
of Conley \cite{Co:1}. We claim that
$$
R(\Phi) = p^{-1}(R_\xi(\phi)) \cup N_0 \cup N_1.
$$
Clearly, $N_0 \cup N_1 \subset R(\Phi)$. By the first part of the proof, we
know that $p^{-1}(R_\xi(\phi))$ is disjoint from and thus has
positive distance to $N_0 \cup N_1$, so that $\rho$ is bounded
below by a positive constant there. Arguing as in the
first part of the proof, one shows that given $n\in N \setminus (N_0
\cup N_1)$, for all sufficiently small $\delta>0$ (depending on $n$)
and all $T\geq 1$ a closed $(\delta,T)$-chain for the flow $\Phi$
based at $n$ cannot enter a sufficiently small neighborhood of
$N_0 \cup N_1$, and so it projects to a closed $(\delta,T')$-chain
based at $x=p(n)$ with vanishing $\lambda$-integral. Similarly, the lift
of any closed $(\delta,T)$-chain based at some $x\in R_\xi(\phi)$ for
sufficiently small $\delta>0$ with vanishing $\lambda$-integral is a
closed $(\delta,T)$-chain based at $n=p^{-1}(x)\in N$ for
$\Phi$. Together these observations prove that $R(\Phi) \setminus (N_0
\cup N_1) = p^{-1}(R_\xi(\phi))$. 

According to the smooth version of Conley's theorem \cite{Co:1}, there
exists a smooth function $L_N:N \lra \bbr$ whose critical set
coincides with $R(\Phi)$ and such that $dL_N(V)<0$ otherwise. By
shifting and rescaling, we can achieve that $L_N(N_1)=1$ and
$L_N(N_0)=0$. The corresponding function $\wt L_N=l_N \circ p^{-1}:M
\lra S^1$ is continuous, and even smooth outside $L^{-1}(l)=p(\del
N)$. It can be smoothed to a function $\wt L:M \lra S^1$ with $d \wt
L(V)<0$ near $L^{-1}(l)$. The differential $d\wt L$ is the required
Lyapunov 1-form for $(V,R_\xi(\phi))$ representing the class $\xi$.
\end{proof}
\subsubsection*{Proof of Proposition~\ref{prop:lyap_exist_r_xi}}

First suppose that there exists a Lyapunov 1-form $\lambda$ for
$(V,R_\xi)$ representing $\xi\in H^1(M;\bbr)$. 

Let $\mu$ be a finite positive invariant measure with $\mu(M\setminus
R_\xi)>0$. It is well known that the support of $\mu$ is contained in
$R$.  Denote by $\mu'$ the restriction of $\mu$ to $M\setminus R_\xi$. 
By Proposition~\ref{prop:flow_cycles}, $\mu'$ gives rise to a
structure cycle $T$ relative to every neighborhood $U$ of $R_\xi$. Choose
$U$ so small that $\xi\in H_U$ and $U$ is disjoint from some
neighborhood of $\supp \mu'$. Then pick a representative $\alpha$ of
$\xi$ which vanishes on $U$ and agrees with $\lambda$ on $\supp
\mu'$. In view of Theorem~\ref{thm:main1} we then have
$$
0> \ct(\xi) = \int_M \alpha(V) d\,\mu' = \int_M \alpha(V) d\mu = A_\mu(\xi).
$$

Conversely, suppose there does not exist a Lyapunov 1-form for
$(V,R_\xi)$ representing the class $\xi$. In the terminology of
\S~\ref{cone-struct}, this means that there is no relatively
transversal form in the cohomology class $-\xi$. Choose a neighborhood
$U$ of $R_\xi$ such that $\ol U$ is disjoint from $R \setminus R_\xi$ and
$\xi$ vanishes on some neighborhood of $\ol U$. Let $T$ be a structure
cycle relative to $U$ with $T_{|\co_U}> 0$ and $T_{|\ca_{U,-\xi}}\leq
0$, which exists by Remark~\ref{rem:small_nbhd} if $U$ is sufficiently
small. By Proposition~\ref{prop:rep_measure} the current $T$ is of the 
form 
$$
T(\alpha) = \int_{M\setminus U} \alpha(V)\, d\mu_T
$$
for some positive Borel measure $\mu_T$. Lemma~\ref{lem:support} now
asserts that $\supp T \subset R\cap(M \setminus U)$. Moreover $\ol U$
is disjoint from $R \setminus R_\xi$, so we conclude from
Lemma~\ref{lem:struct}(ii) that $\del T= 0$, i.e. the measure $\mu_T$
is in fact invariant. Choosing $\alpha\in \ca_{U,-\xi}$ with
$\alpha\equiv 0$ on some neighborhood of $\ol U$, we find that 
$$
A_{\mu_T}(\xi)= - T(\alpha) \geq 0.
$$
This completes the proof of the proposition. \hfill $\Box$

%
%

\section*{Appendix: Proof of Lemma~\ref{lem:approx}}

We briefly recall the setup. $Z$ is an isolated invariant set, $B$ is
an isolating block for it and $Z_0$ is some connected component of
$Z$. For each component $Z'$ of $Z$, we introduced the set $A(Z')$
consisting of those point such that at least one half orbit stays in
$B$ and the corresponding limit set is contained in $Z'$. We denote
the union of these sets by $A(Z)$. Let $\{Z_\iota\}_{\iota\in I}$ be
some indexing of the components of $Z$ different from $Z_0$ ($I$ may
be uncountable).  Finally, denote by $\del_0 B$ the closure of the
complement of $\del_\pm B$ in $\del B$.

The idea now is to construct two sequences of functions $g_n:M \lra
[0,1]$ and $h_n:M \lra [0,1]$ with the following properties:
\begin{enumerate}[(1)]
\item $dg_n(V)=dh_n(V)=0$ on $B$ and $dg_n=dh_n=0$ on some
neighborhood of $Z$.
\item $h_n\equiv 1$ near $A(Z_0)$ and $h_n \lra 0$ pointwise on $B
  \setminus A(Z)$. 
\item $g_n \equiv 1$ near $A(Z_0)$ and $g_n \equiv 0$ near
  $A(Z_\iota)$ whenever $d(A(Z_\iota),A(Z_0))\geq \frac 1 n$, near
  $\del_0 B$ and outside the $\frac 1 n$-neighborhood of $B$.
\end{enumerate}
Then $f_n:=g_n \cdot h_n$ will converge to the characteristic function
of $A(Z_0)$ as required.

Note that given any component $Z_\iota$ of $Z$ different from $Z_0$,
there is a separation $A(Z)=U_\iota \sqcup W_\iota$ into disjoint
closed and open sets with $A(Z_\iota) \subset U_\iota$ and
$A(Z_0)\subset W_\iota$. To construct $g_n$, choose finitely many
indices $\iota_l(n)\in I$ such that the sets $U_{\iota_l}$ cover 
$C_n:=\cup \{A(Z_\iota) \,|\, d(A(Z_\iota),A(Z_0))\geq \frac 1 n\}$.
Then $U_n :=\cup U_{\iota_l(n)}$ and $W_n := \cap W_{\iota_l(n)}$ are
disjoint closed and open sets separating $C_n$ from $A(Z_0)$. 
Set $u_n= (U_n \cup \del_0 B) \cap \del_+ B$ and $w_n =W_n \cap \del_+
B$. Using a smoothing of the distance function from $u_n$, it is now
easy to construct a smooth collared hypersurface $h_n\subset S_+$
separating $w_n$ from $u_n$. From this one first constructs a smooth
function on $S_+$ which equals $1$ near $w_n$ and $0$ near
$u_n$. This can then be extended to a smooth function constant on flow
lines in $B$ which is 0 near $\del_0 B$, and so it can be extended
to a smooth function $g_n$ vanishing outside the $\frac 1 n$-neighborhood of
$B$ which has the properties listed above. 

To construct the functions $h_n$, choose a countable covering of
$\del_+ B \setminus (A(Z) \cup \del_0 B)$ by open disks $D_i$ such
that the disks $\frac 1 2 D_i$ still cover this set. Then first
construct smooth functions on $S_+$ which equal $0$ on $\frac
1 2 D_i$ and equal $1$ outside $\frac 3 4 D_i$ and extend them constantly
along flow lines in $B$  and then to all of $M$ to smooth functions
$\chi_i:M \lra [0,1]$ which equal $1$ outside the $\frac 1
i$-neighborhood of the forward image of $\frac 3 4 D_i$ under the flow
in $B$. Finally define $h_n:= \prod_{i=1}^n \chi_i$. $\Box$

\end{document}